\documentclass{amsart}
\usepackage{amssymb}
\begin{document}
\textheight7in
\begin{abstract} We examine the lower central series of the so-called Fesenko groups $T=T(r)$. These are a certain class of closed subgroups
 of the Nottingham Group.
It is known that all such $T$ are hereditarily just infinite for $p>2$. Here we establish that $T$ has finite width, adding 
to the list of 
known examples. We also prove that $T$ has infinite obliquity. 
\end{abstract}
\title {The Fesenko Groups have finite width}
\author {Cornelius Griffin \\ The Mathematical Institute \\ 24--29 St. Giles' \\ Oxford \\ OX1 3LB}
\thanks{The author thanks the GTEM EU network for their financial support
and the University of Heidelberg for their hospitality while this work was undertaken.
Also thanks to Ivan Fesenko for his helpful comments.}
\address {Mathematical Institute \\ 24--29 St. Giles' \\ University of Oxford \\
Oxford OX1 3LB}
\email {griffinc@maths.ox.ac.uk} 
\date {24th September 2003}
\maketitle

\newtheorem{defn}{Definition}[section]
\newtheorem{thm}{Theorem}[section]
\newtheorem{recipe}{Recipe}[section]
\newtheorem{lemma}{Lemma}[section]
\newtheorem{prop}{Proposition}[section]
\section{Introduction}
\pagestyle{myheadings}
\markboth{Cornelius Griffin}{The Fesenko Groups have finite width}

The Fesenko groups $T:=T(r)$ are defined as follows: let $p$ be a prime number, let $q=p^r$ and define
\[
T  =\left\{ t+\sum_{k\geq 1}a_{qk+1}t^{qk+1}:  a_{qk+1} \in \mathbb F_p\right\} 
\]
with group operation the substitution of one formal power series into another.
These groups are closed subgroups of the Nottingham Group
\[
\mathcal J = \left\{ t+\sum_{k\geq 1}a_kt^{k+1} : a_k \in \mathbb F_p\right\}
\]
which received attention when it was shown in \cite{F} that $T$ form a class of
hereditarily just infinite subgroups of $\mathcal J.$
A pro-$p$ group $G$ is hereditarily just infinite if every open subgroup of $G$ has no infinite quotients.
This is the equivalent condition in pro-$p$ group theory to that of simplicity in finite group theory.
Thus it is natural to hope for a classification of the hereditarily just infinite pro-$p$ groups in a similar manner 
as has been done for the finite simple groups. However it is shown in \cite{KLP} that this is likely to be a vain hope -- there are
too many groups for us to attempt such a classification. 
It is proposed there to add an extra condition -- that of finite width -- and attempt to classify the groups arising now.

\begin{defn} Let $G$ be a pro-$p$ group and write $\gamma_n(G)$ for the $n$-th term of the lower central series of $G.$ Then $G$
  has finite width if there exists
 a constant $C$ independent of $n$ so that 
\[
|\gamma_n(G):\gamma_{n+1}(G)| \leq C \quad \text {for all } n\in \mathbb N.
\]
\end{defn}

Known examples of hereditarily just infinite pro-$p$ groups with finite width include $\mathcal J$,
 some other classes of subgroups of $\mathcal J$ studied by Barnea and Klopsch
\cite{BK,K}, and more traditional examples coming from the theory of the $p$-adic analytic groups. Recently Erchov has uncovered further
classes of subgroups of $\mathcal J$ with the same properties that he uses to prove that $\mathcal J$ is finitely presented \cite{E1,E2}. 

The aim of this article is to prove the following

\begin{thm} The group $T=T(r)$ has finite width for all primes greater than $2.$
\end{thm}

The proof of this result entails realising elements of $T$ as commutators. The calculations are long and involved and it will be helpful to
 introduce the following notation, some
of which we adapt from \cite{BK}.

\begin{defn} \enumerate
\item Let $u(t) = t+ u_{qk+1}t^{qk+1}+\dots $ with $u_{qk+1}\in \mathbb F_p^{\ast}$. Then we write $\delta (u) =k$ and call $\delta (u)$
 the depth of $u;$
\item write $\Gamma_n (T)$ for $\{ \delta (u): u \in \gamma_n (T)\};$
\item write $\Gamma_n^0(T)$ for $\{ k \in \Gamma_n (T): k \equiv 0 \mod p\};$
\item write $\Gamma_n^{-1}(T)$ for $\{ k \in \Gamma_n (T): k \not\equiv 0 \mod p\}$.
\endenumerate
\end{defn}
So notice that we have 
\[
\Gamma_n(T)=\Gamma_n^0 \cup \Gamma_n^{-1}.
\] 
We also define
\[
\Delta_n^{\star}:= \Gamma_n^{\star}\backslash\Gamma_{n+1}^{\star} \quad\text {for } \star \in \{ \ ,0,-1\}.
\]
Then a set of coset representatives for ${\gamma_n(T)}/{\gamma_{n+1}(T)}$ is
\[
\{t+t^{qk+1}: k\in\Delta_n\}
\]
and so $T$ has finite width iff there is a universal bound on $|\Delta_n|.$
Finally note that we will write
$[v,u]=v^{-1}u^{-1}vu.$

\vskip 0.5in
\section{The Calculations}
\vskip 0.25in

To prove that $T$ has finite width it is necessary for us to produce elements of various depths as commutators. We will proceed in two
 steps: firstly we will make note of the
various ways in which these commutators can be generated, and secondly we will piece together these commutators to ensure that 
$|\Gamma_n(T)|$ grows sufficiently slowly for the
results to hold. 

To generate elements as commutators we will use a series of recipes that provide us with the data we need. 
These recipes are taken from \cite{F} where they were first established. They were confirmed by differing methods in \cite{G}.

\begin{recipe} Let $i>j$ and let $(i,p)=1.$ Let $u=t+u_{qi+1}t^{qi+1}, v= t+v_{qk+1}t^{qk+1}+\dots$
 with $u_{qi+1}v_{qj+1}\not=0.$ Then
\[
[u,v]=t+iu_{qi+1}v_{qj+1}t^{q(qj+i)+1}+\dots
\]
and so $[T_i,T_j]\leq T_{qj+i}.$ 
\label{rec1}
\end{recipe}

\begin{recipe} Let $i=i_0p^{n(i)}$ where $p$ does not divide $i_0.$ Let $u=t+u_{qi+1}t^{qi+1}, v= t+v_{qj+1}t^{qj+1}+\dots$ with
$u_{qi+1}v_{qj+1}\not=0$ and suppose that
\[
i>j(\frac {qp^{n(i)}-1}{q-1}).
\] 
Then
\[
[u,v]=t+i_0u_{qi+1}v_{qj+1}t^{q(qp^{n(i)}j+i)+1}+\dots .
\]
\label{rec2}
\end{recipe}  

\begin{recipe} Let $i>j\geq q^2+q$ be so that $i,i-j$ are coprime to $p.$ Then there exists elements $u(l),v(l)$ with 
$\delta (u(l)) \geq i, j-q\leq\delta (v(l)) \leq j$ 
and
with $p$  not dividing $ \delta (u(l)) (\delta (u(l))-\delta (v(l)))$ so that 
\[
t+t^{q(qi+qj)+1}+\dots =\prod_{l=1}^m [u(l),v(l)].
\]
More generally in the same way we may realise the element
\[
t+t^{1+q(qj+p^si)}+\dots
\]
for all $s$ so that $p^s\leq q.$ 
\label{rec3} 
\end{recipe}

 {{\textbf Remark:}} In \cite{F} where Recipe \ref{rec3} is developed Fesenko stipulates that we must also have that $p$ does not divide $j.$
  However the result holds
irregardless of the $p$-adic value of $j$. In any case one can note that this will only impact the arguments that follow for $p=3;$
 for notice that when $p=3$ then if $i,j,i-j 
\not\equiv 0 \mod p$ it follows that $i+j\equiv 0 \mod p$ but for all primes greater than $3$ we can always solve the equation 
$i+j \equiv A \mod p$ with $i,j,i-j \not\equiv 0 \mod p$ for
any $A \in \mathbb N.$ 

\
\newline
It will be useful to consider first the case $p=q$ as this case is the least messy yet still contains all the salient points of the argument
 in the general case. 
\newline
\
\newline
So to begin, information about the elements of depth not divisible by $p$ is contained in the following

\begin{lemma} We have $|\Delta_1^{-1}|=p$ and $|\Delta_n^{-1}| =p-1$ for all $n>1.$ \label{lemma2.1}
\end{lemma}
\proof We use Recipe \ref{rec1}. Notice immediately that we must have
$\gamma_2(T) \leq T_{p+2}$ from which it follows that $ |\Delta_1^{-1}|\geq p.$ Furthermore Recipe 2.1 with $j=1$ and $i>j$ coprime to $p$
 tells us that
\[
[t+t^{1+pi},t+t^{1+pj}]=t+it^{1+p(p+i)}+\dots
\]
and so
\[
\{\lambda : (\lambda ,p)=1, \lambda \geq p+2 \} \subseteq \Gamma_2^{-1}
\]
\[
\implies \Delta_1^{-1}(T) \subseteq \{ \lambda : (\lambda ,p)=1, 1\leq \lambda \leq p+1\}
\]
and so $ |\Delta_1^{-1}|\leq p.$ Thus we have equality.

For $n\geq 1$ we proceed by induction. Suppose that it is known 
\[
\Gamma_n^{-1} =\{(n-1)p+\lambda :(p,\lambda )=1, \lambda \geq 2\}.
\]
Note we have seen this is true for $n=2.$ If it is true for $n$ then for all $x\in \Gamma_n^{-1}$ we have 
\[
[t+t^{1+px},t+t^{1+p}]=t+xt^{1+p(p+x)} +\dots
\]
and so if $x=(n-1)p+\lambda,$ then $np+\lambda \in \Gamma_{n+1}^{-1}.$
The result follows.
\endproof

As we see it is a straightforward task to generate elements of depth prime to $p.$ It is in general more difficult to generate elements
 of depth divisible by $p.$ To demonstrate
how things work we will calculate $\Gamma_2^0$ explicitly.

First notice: for $i>j\geq p^2+p$ so that $i,i-j$ are coprime to $p,$ by Recipe 2.3 we can see that
\[
t+t^{1+p(p(i+j))}+\dots \in \gamma_2(T)
\]
from which it follows that 
\begin{equation}
\{\gamma p: \gamma\geq 2p^2+2p+1\} \subseteq \Gamma_2^0. \label{eq1} 
\end{equation}
Now Recipe 2.2 tells us that given any $\alpha \geq 2$ not divisible by $p$ we have
\[
[t+t^{1+p(\alpha p)},t+t^{1+p}]=t+\alpha t^{1+p(p^2+\alpha p)}+\dots \in \gamma_2(T)
\]
and so
\begin{equation}
\{(\alpha + p)p: \alpha\geq 2, p\not | \alpha\} \subseteq \Gamma_2^0. \label{eq2} 
\end{equation} 
Recipe 2.2 also tells us that given any $\beta \geq 2p$ so that $p$ is the highest power of $p$ dividing $\beta$ we have
\[
[t+t^{1+p(\beta p)},t+t^{1+p}]=t+\beta_0 t^{1+p(p^3+\beta p)}+\dots \in \gamma_2(T)
\]
and so
\begin{equation}
\{(\beta + p^2)p: \beta\geq 2p, p|\beta, p^2\not | \beta\} \subseteq \Gamma_2^0. \label{eq3} 
\end{equation} 

So it follows from (\ref{eq1}) that
\[
\Delta_1^0 \subseteq \{ \gamma p: 1\leq\gamma\leq 2p^2+2p \};
\]
from (\ref{eq2}) that
\[
p(p+\alpha) \not\in\Delta_1^0 \text { for all } \alpha\geq 2 \text { so that } p\not |\alpha;
\]
and from (\ref{eq3}) that
\[
p(p\lambda) \not\in\Delta_1^0 \text { for all } \lambda\geq p+2 \text { so that } p\not |\lambda.
\]
Thus 
\[
\Delta_1^0 \subseteq \{p,2p,\dots,p(p+1); p(2p),p(3p),\dots,p(p^2+p);p(2p^2)\}
\]
and so $|\Delta_1^0|\leq 2p+2$ and $|\Delta_1|\leq 3p+2.$
\newline
\ 
\newline
We now make a more general statement about $\Gamma_n^0:$ in order to ease notation we make the following definitions: firstly we define 
$ C:=\{\gamma p: \gamma \geq 2(p^2+p)+1 \},$ and we also set $A_n:= \{\alpha p: \alpha \geq (n-1)p+2, p\not| \alpha \}$ and 
$B_n:= \{ \beta p: \beta \geq p((n-1)p+2), p^2\not| \beta \}.$ Now we are in a position to prove 

\begin{prop}
The sets $\Gamma_n^0$ for $n\geq 2$ can be described as follows:
\enumerate
\item For $n=2$ we have
\[
\Gamma_n^0 = C\cup A_n\cup B_n;
\]
\item for $n=3$ we have
\[
\Gamma_n^0 = C\cup A_n \cup \{p(2p(p+1)))\};
\]
\item for $4 \leq n<p=2$ we have
\[
\Gamma_n^0 =C\cup A_n:
\]
\item for $n\geq p+2$ we have
\[
\Gamma_n^0 = \{\gamma p: \gamma\geq p^2+p+np+2\}.
\]
\endenumerate
\label{prop2.1}
\end{prop}

\proof
The case $n=2$ follows exactly as the case $n=1$ detailed above. The case $n=3$ is similar: notice that $p(2p(p+1))$ comes from applying
 Recipe 2.2 with $i=p(p(p+2))$ and $j=1.$
 For $4\leq n \leq p+2$ notice that the smallest depth of an element not divisible by $p$ in $\Gamma_{n-1}^{-1}$
is $(n-2)p+2<p^2+p+1$ and so Recipe 2.2 still gives elements not generated by Recipe 2.3. However for $n\geq p+2$ we have that
 $(n-2)p+2\geq p^2+p+1$ and so Recipe 2.2 becomes
redundant.

The result follows.

\endproof

It is now straightforward to deduce the main result in the case $r=1:$

\begin{thm} For $r=1, T=T(1)$ has finite width. Furthermore $w(T) = 3p = |\gamma_1(T):\gamma_2(T)|.$ \label{thm2.1}
\end{thm}

\
\newline
{\textbf Remark.} Notice that it is as expected that for all sufficiently large $n$ we should have $|\gamma_n(T):\gamma_{n+1}(T)| =p^{2p}.$
 To generate elements as commutators,
Recipes 2.1 and 2.3 -- which are the main tools -- only use elements of depth coprime to $p$ which we may always commutate with arbitrary
 elements of $T.$ As Recipe 2.1 tells us 
that
at each stage the depth of elements coprime to $p$ increases by $p$ we should expect a similarly uniform increase in $\Gamma_n.$
 Essentially the same thing will occur when we go
on to consider the case $q=p^r.$ The generation of elements of depth coprime to $p$ follows from Recipe 2.1 exactly as in the case $p=q.$
 However 
now we will use Recipe 2.3 to generate elements of depth $q\lambda$ in $\gamma_n(T)$ for all $\lambda$ greater than some $\lambda (n).$
 We will
also show that $\lambda(n+1)=\lambda (n) +q$ for all sufficiently large $n.$ Then we will use Recipe 2.2 to generate sufficient elements of
 depth divisible by $p^s$ for 
$s<r$ to complete the proof.

\vskip 0.5in
\section{The General Case}
\vskip 0.25in

As previously it is easy to gather together the information about elements of depth not divisible by $p:$

\begin{lemma} We have $|\Delta_1^{-1} (T)| = q-p^{r-1}+1$ and $|\Delta_n^{-1}(T)|= q-p^{r-1}$ for all $n>1.$
\end{lemma}

\proof
As in Lemma \ref{lemma2.1} we have
\[
\Gamma_n^{-1} =\{\lambda :(\lambda,p)=1, \lambda\geq (n-1)q+2\}
\]
from the commutator
\[
[t+t^{1+q((n-1)q+\mu)},t+t^{1+q}]=t+\mu t^{1+q(nq+\mu)}+\dots .
\]
Thus
\[
|\Delta_1^{-1}|=|\{\lambda: (\lambda,p)=1, 1\leq\lambda\leq q+1\}| = q+1-p^{r-1}
\]
and similarly for $n>1,$
\[
|\Delta_n^{-1}|= |\{\lambda:(\lambda,p)=1, (n-1)q+2\leq \lambda<nq+2 \}|=q-p^{r-1}.
\]

\endproof
Notice that this agrees with Lemma \ref{lemma2.1} in the case $r=1.$

We now proceed in two stages: firstly we show how to generate elements of depth divisible by $q$, and secondly we show how to generate
 elements of depth divisible by $p^s$ for
$s<r.$ It will be helpful to introduce some more notation. We set
\[
\epsilon_n^s(T):=\{\delta (u):u \in \gamma_n (T), p^s|\delta(u), p^{s+1}\not|\delta (u)\} \text { for all } s<r
\]
and
\[
\epsilon_n^r(T):=\{\delta (u):u \in \gamma_n (T), p^r|\delta(u) \} .
\]
Then all the information we need is contained in the following

\begin{prop}
With the above notation:
\enumerate
\item for all $s<r$ and all $n$
\[
\epsilon_n^s(T) =\{p^s\lambda : \lambda \geq (n-1)q+1, p \text { does not divide } \lambda\};
\]
\item for $n< q+2$
\[
\epsilon_n^r(T) =\{\gamma q: \gamma\geq 2q^2+2q+1\};
\]
\item for $n\geq q+2$
\[
\epsilon_n^r(T) =\{\gamma q: \gamma \geq q^2+q+nq+2\}.
\]
\endenumerate
\label{prop3.1}
\end{prop}

\proof
To get (1) notice that for $i_0$ coprime to $p,$
\[
[t+t^{1+q(i_0p^s)}, t+t^{1+q}]= t+t^{1+q(qp^s+i_0p^s)}+\dots 
\]
from which (1) follows easily by induction. i.e., assume the result is known for $k \leq n$ and then use the relationship
\[
[t+t^{1+q(p^s((n-1)q+\mu)},t+t^{1+q}] = t+\mu t^{1+q(p^s(nq +\mu))}+\dots 
\]
for $\mu$ not divisible by $p.$

(2) and (3) follow as in Proposition \ref{prop2.1} from Recipe 2.3. Notice that 
\[
(n-1)q+2\geq q^2+q+1
\]
\[
\iff n\geq q+2-\frac 1q >q+1.
\]

\endproof

From here it follows automatically that we may deduce
\begin{thm} The groups $T=T(r)$ have finite width for $r\geq 2$ and so combined with Theorem 2.1 for $r\geq 1.$
\end{thm}

\proof
Notice that 
\[
|\Delta_n(T)| = |\Delta_n^{-1}(T)| +\sum_{s=1}^r |\epsilon_n^s(T)\backslash \epsilon_{n+1}^s(T)|
\]
and so one simply has to combine the results we have proved above in Lemma 3.1 and Proposition 3.1.

\endproof

We now go on to show that $T$ has infinite obliquity.
\vskip 0.5in
\section{The Obliquity of $T$} 
\vskip 0.25in

We recall the definition of obliquity in relation to a pro-$p$ group $G$. Write 
\[
\mu_n(G) = \gamma_{n+1}(G)\cap \{N\triangleleft G: N\not< \gamma_{n+1}(G)\}
\] 
and then the obliquity of $G$ is defined to be
\[
o(G) =\sup_n \log_p |\gamma_{n+1}(G):\mu_n(G)| \in \mathbb N \cup\infty.
\]
So notice that a group $G$ has obliquity $0$ iff and only if given any normal subgroup $N$ of $G$ there exists some integer $n$ so that 
$\gamma_{n+1}(G)\leq N \leq\gamma_n(G).$
More generally $G$ has finite obliquity $t$ if and only if given a open normal subgroup $N$ of $G$ then there exists $n$ so that 
$\gamma_{n+t+1}(G) \leq N \leq \gamma_n(G)$
\cite{CC}. We show now that $T$ is as far removed from this as possible:

\begin{prop} $T$ has infinite obliquity.
\end{prop}

\proof
It is necessary to produce, for every $n \in \mathbb N$ a normal subgroup $N(n)$ of $T$ so that $N\not< \gamma_{n+1}(T)$ and so that
\[
|\gamma_{n+1}(T):N\cap\gamma_{n+1}(T)| \rightarrow \infty \text { as } n\rightarrow \infty.
\]
We will use the following two facts:
\enumerate
\item $\inf \{\delta (\gamma_n(T))\} = (n-1)q+2;$
\item  $\inf \{\delta (u):u\in \gamma_n(T), p\text { divides }\delta (u)\} = p((n-1)q+1).$
\endenumerate

Given this, set 
\[
H_n = T_{p((n-1)q)+1}.
\] 
So as $t+t^{q(p((n-1)q+1))+1}\in \gamma_{n+1}(T), \not\in H_n$ it follows from (2) that $H_n$ is not contained in $\gamma_{n+1}(T)$ and
 from (1) and Lemma 3.1 that
\[
\lim_{n\rightarrow \infty} |\gamma_{n+1}(T):\gamma_{n+1}(T)\cap H_n| =\infty.
\]
The result follows. 

\endproof

\end{document}